\documentclass[11pt]{elsarticle} % use larger type; default would be 10pt
\usepackage[utf8]{inputenc} % set input encoding (not needed with XeLaTeX)

\usepackage{amsmath, amssymb, amsfonts, float, color}
\usepackage{subfig,epstopdf}

\newtheorem{theorem}{Theorem}[section]

\newtheorem{proposition}[theorem]{Proposition}

\newproof{pf}{Proof}

\usepackage{booktabs} % for much better looking tables
\usepackage{array} % for better arrays (eg matrices) in maths
\usepackage{paralist} % very flexible & customisable lists (eg. enumerate/itemize, etc.)
\usepackage{verbatim} % adds environment for commenting out blocks of text & for better verbatim
\usepackage{subfig} % make it possible to include more than one captioned figure/table in a single float
\DeclareMathOperator*{\argmin}{arg\,min}
\DeclareMathOperator{\vspan}{span}

\begin{document}

\begin{frontmatter}
\title{State Following (StaF) Kernel Functions for Function Approximation\tnoteref{t1}}

\tnotetext[t1]{A preliminary version of this work appeared in the proceedings of the 2015 American Control Conference (ACC)  \cite{Rosenfeld.Kamalapurkar.ea2015}}

\author[myu]{Joel A. Rosenfeld\corref{cor1}}
\ead{joelar@ufl.edu}

\author[myu]{Rushikesh Kamalapurkar}

\author[myu]{Warren E. Dixon\thanks{test}}

\address[myu]{Department of Mechanical and Aerospace Engineering, MAE-B, University of Florida, Gainesville, FL}

\cortext[cor1]{Corresponding author}

\begin{abstract}
A function approximation method is developed that aims to approximate a function in a small neighborhood of a state that travels within a compact set. The development is based on the theory of universal reproducing kernel Hilbert spaces over the $n$-dimensional Euclidean space. Several theorems are introduced that support the development of this State Following (StaF) method. In particular, it is shown that there is a bound on the number of kernel functions required for the maintenance of an accurate function approximation as a state moves through a compact set. Additionally, a weight update law, based on gradient descent, is introduced where arbitrarily close accuracy can be achieved provided the weight update law is iterated at a sufficient frequency, as detailed in Theorem \ref{thm:gradient-chase}.

To illustrate the advantage, the impact of the StaF method is that for some applications the number of basis functions can be reduced. The StaF method is applied to an adaptive dynamic programming (ADP) application to demonstrate that stability is maintained with a reduced number of basis functions.

Simulation results demonstrate the utility of the StaF methodology for the maintenance of accurate function approximation as well as solving an infinite horizon optimal regulation problem through ADP. The results of the simulation indicate that fewer basis functions are required to guarantee stability and approximate optimality than are required when a global approximation approach is used.
%\footnote{A preliminary version of this work appeared in the proceedings of the 2015 American Control Conference (ACC)  \cite{Rosenfeld.Kamalapurkar.ea2015}}
\end{abstract}

\end{frontmatter}
%\maketitle
\section{Introduction}
% Dixon wants to move all ADP discussion to a section right before the ``conclusion'' section
% This includes HJB optimal control problems etc.
% The intro should lead with moving local approximations following the state of a system.
% Mention ADP as a motivation, but don't discuss until much later.
% In the ADP Section - Mention how the theory is constrained to approximations that achieve global 
% approximations, and in implementation only polynomials have been used, since it is difficult to 
% achieve weight convergence for the large number of basis functions required by other methods.
% A severe limitation in the field is ... simulations .... blah blah...
%
% Intro should mention that this is a method to reduce the number of basis functions, that differs 
% from other literature in that it focuses on the establishment and maintenance of a local 
% approximation.
Often in the theory of approximation, an accurate estimation of a function over a large compact set is sought \cite{Christmann.Christmann.ea2010, Micchelli.Micchelli.ea2007, Park.Park.ea1991}. It is well known that the larger the compact set, a correspondingly larger number of basis functions are required to achieve an accurate function approximation. There is a large body of literature concerned with methods for the reduction of the number of basis functions required to achieve such an approximation (c.f. \cite{C.Prudhomme.C.Prudhomme.ea2001, Balmes.Balmes.ea1996, AhmedK.Noor.AhmedK.Noor.ea1980}).

In many control applications, function approximation is used to generate a stabilizing controller of a state in a dynamical system. For instance, in adaptive dynamic programming (ADP), an approximation of the optimal value function is leveraged to produce an approximate optimal controller \cite{Al-Tamimi2008, Bhasin.Kamalapurkar.ea2013a, Dierks2009, Lewis.Vrabie2009, Mehta.Meyn2009, Padhi2006, Vamvoudakis2010, Zhang.Cui.ea2013, Zhang.Cui.ea2011, Zhang.Liu.ea2013}. Traditionally, the approximation is sought over a large compact set, and requires many basis functions. The computational resources required to tune the weights of the basis functions renders real-time implementation of controllers based on ADP methods infeasible.

Motivated by problems in control theory, this paper introduces an approximation methodology that aims to establish and maintain an accurate approximation of a function in a neighborhood of a moving state in a dynamical system. The method, deemed the state following (StaF) method, reduces the number of basis functions required to achieve an accurate approximation by focusing on the approximation of a function over a small neighborhood by linear combinations of time and state varying basis functions. Therefore, even in cases where processing power of on-board CPUs is limited, an accurate approximation of a function can be maintained.

%The aim of this paper is the development of an efficient scheme for the approximation of continuous functions via state and time varying basis functions that maintain the approximation of a function in a local neighborhood of the state, deemed the state following (StaF) method. The method developed in this paper is presented as a general strategy for function approximation, and can be implemented in contexts outside of optimal control. An example of a successful application of this method to ADP is presented in Appendix \ref{sec:adp-simulation}.

The particular basis functions that will be employed throughout this paper are derived from kernel functions corresponding to RKHSs. In particular, the centers are selected to be continuous functions of the state variable bounded by a predetermined value. That is, given a compact set $D \subset \mathbb{R}^n$, $\epsilon > 0$, $r > 0$ and $M \in \mathbb{N}$, $c_i(x) =  x + d_i(x)$ where $d_i:\mathbb{R}^n \to \mathbb{R}^n$ is continuously differentiable and $\sup_{x \in D} \| d_i(x) \| < r$ for $i=1,...,M$. The parameterization of a function $V:D\to \mathbb{R}$ in terms of StaF kernel functions is given by $$\hat V(y;x(t),t) = \sum_{i=1}^M w_i(t) K(y,c_i(x(t)))$$ where $w_i(t)$ is a weight signal chosen to satisfy $$\limsup_{t \to \infty} E_r(x(t),t)< \epsilon$$ where $E_r$ is a measure of the accuracy of an approximation in a neighborhood of $x(t)$, such as that of the supremum norm: $$E_r(x(t),t)=\sup_{y \in \overline{N_r(x(t))}} \left| V(y) - \hat V(y;x(t),t) \right|.$$
 
The goal of the StaF method is to establish and maintain an approximation of a function in a neighborhood of the state. The justification for this approach stems from the observation that an optimal controller only requires the value of the estimation of the optimal value function to be accurate at the current system state. Thus, when computational resources are limited, computational efforts should be focused on improving the accuracy of approximations near the system state.

The advantage of using RKHSs for the purpose of local approximations is twofold. RKHSs have been found to be effective for nonlinear function approximation \cite{DeVore.DeVore.ea1998}, and the use of RKHS can enable accurate estimations of a wide array of nonlinear functions. Also, the ideal weights corresponding to the Hilbert space norm provided by RKHSs change smoothly with respect to smooth changes in the centers, as demonstrated in Theorem \ref{thm:continuous-ideal-weight}, which allows the execution of weight update laws to achieve and maintain an accurate approximation. The ideal weights in the context of the StaF approximation method become a continuous function of the state and are investigated in Section \ref{sec:ideal-weight-function}.

Previous efforts in the literature have performed nonlinear approximation through the adjustment of the centers of radial basis functions (c.f. \cite{Gaggero.Gaggero.ea2013, Gaggero.Gaggero.ea2014, Zoppoli.Zoppoli.ea2002}) as a means to determine the optimal centers for global approximation. These efforts are more applicable when off-line techniques can be used due to computational demands. For other applications where computational resources are limited, global approximations may not be feasible (especially as the dimension of the problem grows), nor is the optimal selection of parameters.
 
This paper lays the foundation for the establishment and maintenance of a real-time moving local approximation of a continuous function. Section \ref{sec:staf-problem} of this paper frames the particular approximation problem of the StaF method. Section \ref{sec:feasibility} demonstrates accurate approximation with a fixed number of moving basis functions. Section \ref{sec:exp-kern} demonstrates an explicit bound on the number of required StaF basis functions for the case of the exponential kernel function. The ideal weight function arising from the StaF method is introduced and discussed in Section \ref{sec:ideal-weight-function}, where the existence and smoothness of the ideal weight function is established. Section \ref{sec:gradient-chase} provides a proof of concept demonstrating the existence of weight update laws to maintain an accurate approximation of a function in a local neighborhood, ultimately establishing a uniform ultimate bounded result. The remaining sections demonstrate the developed method through numerical experiments and discussions of applications. Specifically, Section \ref{sec:gradient-chase-simulation} gives the results of a ``gradient chase'' algorithm. In Section \ref{sec:ADP-Outline}, the utility of StaF methods are demonstrated in an ADP application.
 
\section{The StaF Problem Statement}
\label{sec:staf-problem}
 
Given a continuous function $V : \mathbb{R}^n \to \mathbb{R}$, $r > 0$, an arbitrarily small $\epsilon > 0$, and a dynamical system $\dot x = f(x,u)$ (where $f$ is sufficiently regular for the system to be well defined), the goal of the StaF approximation method is to select state and time varying basis functions $\sigma_i:\mathbb{R}^n \times \mathbb{R}^n \times \mathbb{R} \to \mathbb{R}$ for $i=1,2,...,M$ and weight signals $w_i : \mathbb{R}_+\to\mathbb{R}$ for $i=1,2,...,M$  such that 
\begin{equation}\label{eq:staf-problem}\limsup_{t\to\infty} \sup_{y \in \overline{N_r(x(t))}} \left| V(y) - \sum_{i=1}^M w_i(t) \sigma_i(y;x(t),t)\right| < \epsilon.\end{equation}
In other words, the StaF approximation method aims to achieve an arbitrarily small steady state error of order $\epsilon$ in a closed neighborhood of the state, $\overline{N_r(x(t))} = \{ y \in \mathbb{R}^n : \| x(t)-y \|_2 \le r \}$.

Central problems to the StaF method include determining the basis functions and the weight signals. When reproducing kernel Hilbert spaces are used for basis functions, \eqref{eq:staf-problem} can be relaxed to where the supremum norm is replaced with the Hilbert space norm. Since the Hilbert space norm of a RKHS dominates the supremum norm, \eqref{eq:staf-problem} with the supremum norm is simultaneously satisfied. Moreover, when using a RKHS, the basis functions can be selected to correspond to centers placed in a moving neighborhood of the state. In particular, given a kernel function $K:\mathbb{R}^n \times \mathbb{R}^n \to \mathbb{R}$ corresponding to a (universal) RKHS, $H$, and continuous center functions $c_i:\mathbb{R}^n \to \mathbb{R}^n$ for which $d_i(x):=c_i(x)-x$ is bounded by $r$, then the StaF problem becomes the determination of weight signals $w_i:\mathbb{R}_+ \to \mathbb{R}$ for $i=1,...,M$ such that:

\begin{equation}\label{eq:staf-problem-rkhs}\limsup_{t\to\infty} \left\| V(\cdot) - \sum_{i=1}^M w_i(t) K(\cdot,c_i(x(t)))\right\|_{r,x(t)} < \epsilon\end{equation} where $\|\cdot\|_{r,x(t)}$ is the norm of the RKHS obtained by restricting functions in $H$ to $N_r(x(t))$.

Since \eqref{eq:staf-problem-rkhs} implies \eqref{eq:staf-problem}, the focus of this paper is to demonstrate the feasibility of satisfying \eqref{eq:staf-problem-rkhs}.  Theorem \ref{thm:feasibility} demonstrates that under a certain continuity assumption a bound on the number of kernel functions necessary  for the maintenance of an approximation throughout a compact set can be determined, and Theorem \ref{thm:continuous-ideal-weight} shows that a collection of continuous ideal weight functions can be determined to satisfy \eqref{eq:staf-problem-rkhs}. Theorem \ref{thm:continuous-ideal-weight} justifies the use of weight update laws for the maintenance of an accurate function approximation, and this is demonstrated by Theorem \ref{thm:gradient-chase} as well as the numerical results contained in Section \ref{sec:gradient-chase-simulation} and \ref{sec:adp-simulation}. 

The choice of RKHS for Section \ref{sec:gradient-chase-simulation} is that which corresponds to the exponential kernel $K(x,y) = \exp(x^Ty)$ where $x, y \in \mathbb{R}^n$ and will be denoted by $F^2(\mathbb{R}^n)$ since it is closely connected to the Bargmann-Fock space \cite{Zhu.Zhu.ea2012}. The RKHS corresponding to the exponential kernel is a universal RKHS \cite{Pinkus.Pinkus.ea2004, Steinwart.Christmann2008}, which means that given any compact set $D \subset \mathbb{R}^n$, $\epsilon >0$ and continuous function $f:D\to \mathbb{R}$, there exists a function $\hat f \in F^2(\mathbb{R}^n)$ for which $\sup_{x \in D} | f(x) - \hat f(x) | < \epsilon$.

\section{Feasibility of the StaF Approximation and the Ideal Weight Functions}
\label{sec:feasibility}

The first theorem concerning the StaF method demonstrates that if the state variable is constrained to a compact subset of $\mathbb{R}^n$, then there is a finite number of StaF basis functions required to establish the accuracy of an approximation.

\begin{theorem}\label{thm:feasibility}Suppose that $K:X\times X \to \mathbb{C}$ is a continuous kernel function corresponding to a RKHS, $H$, over a set $X$ equipped with a metric topology. If $V \in H$, $D$ is a compact set of $X$, $r> 0$, and $\|V\|_{x,r}$ is continuous with respect to $x$, then for all $\epsilon > 0$ there is a $M \in \mathbb{N}$ such that for each $x\in D$ there are centers $c_1,c_2,...,c_M \in N_r(x)$ and weights $w_i \in \mathbb{C}$ such that $$\left\| V(\cdot) - \sum_{i=1}^M w_i K(\cdot,c_i) \right\|_{r,x} < \epsilon.$$\end{theorem}

\begin{pf}
Given $\epsilon > 0$, for each neighborhood $N_r(x)$ with $x\in D$, there exists a finite number of centers $c_1,...,c_M \in N_r(x)$, and weights $w_1,...,w_M \in \mathbb{C}$, such that $$\left\| V(\cdot) - \sum_{i=1}^M w_i K(\cdot, c_i)\right\|_{r,x} < \epsilon.$$ Let $M_{x,\epsilon}$ be the minimum such number. The claim of the proposition is that the set $Q_\epsilon := \{ M_{x,\epsilon} : x \in D\}$ is bounded. Assume by way of contradiction that $Q_\epsilon$ is unbounded, and take a sequence $\{x_n\} \subset D$ such that $M_{x_n,\epsilon}$ is a strictly increasing sequence and $x_n \to x$ in $D$. It is always possible to find such a convergent sequence, since every compact subset of metric space is sequentially compact.
Let $c_1,...,c_{M_{x,\epsilon/2}} \in N_r(x)$ and $w_1,...,w_{M_{x,\epsilon/2}} \in \mathbb{C}$ be centers and weights for which
\begin{equation}\label{eq:epsilon-over-two}E(x):=\left\|V(\cdot) - \sum_{i=1}^{M_{x,\epsilon/2}} w_i K(\cdot, c_i)\right\|_{r,x} < \epsilon/2.\end{equation}
For convenience, let each $c_i \in N_r(x)$ be expressed as $x+d_i$ for $d_i \in N_r(0)$. The function $E(x)$ in \eqref{eq:epsilon-over-two} can be written as
$$\left( \|V\|_{r,x} - 2 Re \left( \sum_{i=1}^{M_{x,\epsilon/2}} w_i V(x+d_i) \right) + \sum_{i,j=1}^{M_{x,\epsilon/2}} w_i \overline{w_j} K(x+d_i, x+d_j) \right)^{1/2}.$$ By the hypothesis, $K$ is continuous with respect to $x$, which implies that $V$ is continuous \cite{Christmann.Christmann.ea2010}, and $\|V\|_{r,x}$ is continuous with respect to $x$. Hence, there exists $\eta > 0$ for which $|E(x)-E(x_n)| < \epsilon/2$ for all $x_n \in N_{\eta}(x)$. Thus $E(x_n) < E(x) + \epsilon/2 < \epsilon$ for sufficiently large $n$. By minimality $M_{x_n,\epsilon} < M_{x,\epsilon/2}$ for sufficiently large $n$. This is a contradiction.\qed\end{pf}

The assumption of the continuity of $\| V \|_{r,x}$ in Theorem \ref{thm:feasibility} is well founded. There are several examples where the assumption is known to hold. For instance, if the RKHS is a space of real entire functions, as it is for the exponential kernel, then $\| V\|_{r,x}$ is not only continuous, but it is constant.

Using a similar argument as that in Theorem \ref{thm:feasibility}, the theorem can be shown to hold when the restricted Hilbert space norm is replace by the supremum norm over $\overline{N_r(x)}$. The proof of the following theorem can be found in the preliminary work for this article in \cite{Rosenfeld.Kamalapurkar.ea2015}.

\begin{proposition}\label{prop:sup-norm-feasibility}Let $D$ be a compact subset of $\mathbb{R}^n$, $V:\mathbb{R}^n \to \mathbb{R}$ be a continuous function, and $K: \mathbb{R}^n \to \mathbb{R}^n \to \mathbb{R}$ be a continuous and universal kernel function. For all $\epsilon, r >0$, there exists $M \in \mathbb{N}$ such that for each $x\in D$, there is a collection of centers $c_1,..., c_M \in N_r(x)$ and weights $w_1,...,w_M \in \mathbb{R}$ such that $\sup_{y \in \overline{N_r(x)}} \left| V(y) - \sum_{i=1}^M K(y,c_i) \right| < \epsilon.$\end{proposition}

\section{Explicit Bound for the Exponential Kernel}
\label{sec:exp-kern}

Theorem \ref{thm:feasibility} establishes a bound on the number of kernel functions required for the maintenance of the accuracy of a moving local approximation. However, the proof does not provide an algorithm to computationally determine the upper bound. Even when the approximation with kernel functions is performed over a fixed compact set, a general bound for the number of collocation nodes required for accurate function approximation under the Hilbert space norm is unknown. Thus, it is desirable to have a computationally determinable upper bound to the number of StaF basis functions required to yield an arbitrarily close approximation. Theorem \ref{thm:exponential-bound} provides a calculable bound on the number of exponential functions required to yield an arbitrarily close approximation with respect to the supremum norm. That is, Theorem \ref{thm:exponential-bound} provides a computable analogue of Theorem \ref{thm:feasibility} and Proposition \ref{prop:sup-norm-feasibility} for a StaF approximation problem of the form $$\limsup_{t \to \infty} \sup_{y \in \overline{N_r(x(t))}} \left| V(y) - \sum_{i=1}^\infty w_i(t) K(y, c_i(x(t)))\right| < \epsilon.$$

While error bounds have been computed for the exponential function with respect to the supremum norm (c.f. \cite{Beylkin.Beylkin.ea2005}), current literature allows the ``frequencies'' or centers of the exponential kernel functions to be unconstrained. The lack of constraints on the centers of the exponential kernel functions means that the existing results cannot be leveraged for the StaF approximation problem. The contribution of Theorem \ref{thm:exponential-bound} is the development of an error bound while constraining the size of the centers.

\begin{theorem}\label{thm:exponential-bound}Let $K:\mathbb{R}^n \times \mathbb{R}^n \to \mathbb{R}$ given by $K(x,y)=\exp\left(x^Ty\right)$ be the exponential kernel function. Let $D \subset \mathbb{R}^n$ be a compact set, $V:D\to \mathbb{R}$ continuous, and $\epsilon, r > 0$. For each $x \in D$, there exists a finite number of centers $c_1,...,c_{M_{x,\epsilon}} \in N_r(x)$ and weights $w_1,w_2,...,w_{M_{x,\epsilon}} \in \mathbb{R}$, such that $$\sup_{y\in \overline{N_r(x)}} \left| V(y) - \sum_{i=1}^{M_{x,\epsilon}} w_i K(y,c_i)\right| < \epsilon.$$ If $p$ is an approximating polynomial that achieves the same accuracy over $\overline{N_r(x)}$ with degree $N_{x,\epsilon}$, then an asymptotically similar bound can be found with $M_{x,\epsilon}$ kernel functions, where $M_{x,\epsilon} < {n+N_{x,\epsilon}+S_{x,\epsilon} \choose N_{x,\epsilon} + S_{x,\epsilon}}$ for some constant $S_{x,\epsilon}$ that is the degree of an approximating polynomial. Moreover, $N_{x,\epsilon}$ and $S_{x,\epsilon}$ can be bounded uniformly over $D$, and thus, so can $M_{x,\epsilon}$.\end{theorem}

\begin{pf}For notational simplicity, the quantity $\| f \|_{D,\infty}$ denotes the supremum norm of a function $f: D \to \mathbb{R}$ over the compact set $D$ throughout the proof of Theorem \ref{thm:exponential-bound}.

First, consider the ball of radius $r$ centered at the origin. The statement of the theorem can be proved by finding an approximation of monomials by a linear combination of exponential kernel functions.

Let $\alpha=(\alpha_{1},\alpha_{2},...,\alpha_{n})$ be a multi-index, and define $|\alpha|=\sum\alpha_{i}$. Note that\footnote{The notation $g_{m}(x)=O(f(m))$ means that for sufficiently large $m$, there is a constant $C$ for which $g_{m}(x)<Cf(m)$ for all $y\in\overline{N_{r}(0)}$.} 
$$
m^{|\alpha|}\prod_{i=1}^{n}\left(\exp\left(y_{i}/m\right)-1\right)^{\alpha_{i}}=y_{1}^{\alpha_{1}}y_{2}^{\alpha_{2}}\cdots y_{n}^{\alpha_{n}}+O\left(\frac{1}{m}\right)
$$
which by the binomial theorem leads to the sum
\begin{multline}
m^{|\alpha|}\sum_{l_{i}\le\alpha_{i},i=1,2,...,n}{\alpha_{1} \choose l_{1}}{\alpha_{2} \choose l_{2}}\cdots{\alpha_{n} \choose l_{n}}(-1)^{|\alpha| - \sum_{i}l_{i}} \exp\left(\sum_{i=1}^{n}y_{i} \left(\frac{l_{i}}{m}\right)\right)\\
=y_{1}^{\alpha_{1}}y_{2}^{\alpha_{2}}\cdots y_{n}^{\alpha_{n}}+O\left(\frac{1}{m}\right).\label{eq:monomialapprox}
\end{multline}

The big-oh constant indicated by $O(1/m)$ can be computed in terms of the derivatives of the exponential function via Taylor's Theorem. The centers corresponding to this approximation are of the form $l_{i}/m$ where $l_{i}$ is a non-negative integer satisfying $l_{i}<\alpha_{i}$. Hence, for $m$ sufficiently large, the centers reside in $N_{r}(0)$. 

To shift the centers so that they reside in $N_{r}(y)$, let $x=(x_{1},x_{2},...,x_{n})^{T}\in\mathbb{R}^{n}$,
and multiply both sides of (\ref{eq:monomialapprox}) by $\exp\left(y^{T}x\right)$ to get 

$$m^{|\alpha|}\sum_{l_{i}\le\alpha_{i},i=1,2,...,n}{\alpha_{1} \choose l_{1}}{\alpha_{2} \choose l_{2}}\cdots{\alpha_{n} \choose l_{n}}(-1)^{|\alpha|-\sum_{i}l_{i}}\exp\left(\sum_{i=1}^{n}y_{i}\left(\frac{l_{i}}{m}+x_{i}\right)\right)$$
$$=e^{y^{T}x}\left(y_{1}^{\alpha_{1}}y_{2}^{\alpha_{2}}\cdots y_{n}^{\alpha_{n}}\right)+O\left(\frac{1}{m}\right).$$

For each multi-index, $\alpha=(\alpha_{1},\alpha_{2},...,\alpha_{n})$, the centers for the approximation of the corresponding monomial are of the form $x_{i}+l_{i}/m$ for $0\le l_{i}\le\alpha_{i}$. Thus, by linear combinations of these kernel functions, a function of the
form $e^{y^{T}x}g(y)$, where $g$ is a multivariate polynomial, can be uniformly approximated by exponential functions over $N_{r}(x)$. Moreover if $g$ is a polynomial of degree $\beta$, then this approximation can be a linear combination of ${n+\beta \choose \beta}$ kernel functions.

Two polynomials, $p_x$ and $q_x$ are selected to approximate $V$ and $e^{-x^Ty}$, respectively, over $N_r(x)$. Since $V$ is a continuous function, it can be approximated with arbitrary accuracy by polynomials. Subsequently, the previous development will be utilized to approximate the polynomials by linear combinations of exponential functions.

Let $\epsilon'>0$ and suppose that $p_{x}$ is polynomial with degree $N_{x,\epsilon'}$ such that $$p_{x}(y)=V(y)+\epsilon_{1}(y)$$ where $|\epsilon_{1}(y)|<\|e^{y^{T}x}\|_{D,\infty}^{-1}\epsilon'/2$ for all $y\in N_{r}(x)$. Let $q_{x}(y)$ be a polynomial in $\mathbb{R}^{n}$ variables of degree $S_{x,\epsilon}$ such that $$q_{x}(y)=e^{-y^{T}x}+\epsilon_{2}(y)$$ where $\epsilon_{2}(y)<\|V\|_{D,\infty}^{-1}\|e^{y^{T}x}\|_{D,\infty}^{-1}\epsilon'/2$ for all $y\in N_{r}(x)$.

The above construction indicates that there is a sequence of linear combinations of exponential kernel functions, $F_{m}(y)$, (with a fixed number of centers inside $N_{r}(x)$) for which 
\begin{align*}F_{m}(y)&=e^{y^{T}x}q_{x}(y)p_{x}(y)+O\left(\frac{1}{m}\right)\\&=e^{y^{T}x}\left(e^{-y^{T}x}+\epsilon_{2}(y)\right)\left(V(y)+\epsilon_{1}(y)\right)+O\left(\frac{1}{m}\right).\end{align*}
After multiplication and an application of the triangle inequality, the following is established: 
$$\left|F_{m}(y)-V(y)\right|<\epsilon'+\left(\frac{\|V\|_{D,\infty}^{-1}\|e^{y^{T}x}\|_{D,\infty}^{-1}}{4}\right)\epsilon'^{2}+O\left(\frac{1}{m}\right)$$
 for all $y\in N_{r}(x)$. The degree of the polynomial $q_{x}$, $S_{x,\epsilon}$, can be uniformly bounded in terms of the modulus of continuity of $e^{y^{T}x}$ over $D$. Similarly, the uniform bound on the degree of $p_{x}$, $N_{x,\epsilon'}$, can be described in terms of the modulus of continuity of $V$ over $D$. The number of centers required for $F_{m}(y)$ is determined by the degree of the polynomial $q\cdot p$ (treating the $x$ terms of $q$ as constant), which is sum of the two polynomial degrees. Finally for $m$ large enough and $\epsilon'$ small enough, $|F_{m}(y)-V(y)|<\epsilon$, and the proof is complete.
\qed\end{pf}

Theorem \ref{thm:exponential-bound} demonstrates an upper bound required for the accurate approximation of a function through the estimation of approximating polynomials. Moreover, the upper bound is a function of the polynomial degrees. For example, for a neighborhood of the origin in $\mathbb{R}$, if $p$ is an approximating polynomial of degree $N$, then the same order of approximation can be achieved by a linear combination of $N+1$ exponential functions. The exponential kernel will be used for simulations in Section \ref{sec:gradient-chase-simulation} and \ref{sec:adp-simulation}.

\section{Existence and Smoothness of the Ideal Weight Function}\label{sec:ideal-weight-function}

Theorem \ref{thm:feasibility} and Proposition \ref{prop:sup-norm-feasibility} establish that given a kernel function, a finite number of centers can be used to yield an arbitrarily accurate estimation of a function, for a set of ideal weights. Theorem \ref{thm:exponential-bound} further establishes that for the exponential kernel function, a calculable number of centers can be determined. However, further investigation is required to understand the characteristics of the ideal weights that correspond to the moving centers. For example, in control applications involving function approximation or system identification, it is assumed that there is a collection of constant ideal weights, and much of the theory is in the demonstration of the convergence of approximate weights to the ideal weights. The subsequent Theorem \ref{thm:continuous-ideal-weight} establishes that ideal weights, which are functions of the state dependent centers, are $m$-times continuously differentiable. This property can then be used to develop weight update laws (e.g., see Section \ref{sec:gradient-chase}).%Since the ideal weights are no longer constant, it is necessary to show that the ideal weights change smoothly as the system progresses. The smooth change in centers will allow the proof of uniform ultimately bounded (UUB) results through the use of weight update laws. One such result will be demonstrated in Section \ref{sec:gradient-chase}, in particular a UUB result is proven in Theorem \ref{thm:gradient-chase}.

Since the ideal weights corresponding to a Hilbert space norm are unique, Theorem \ref{thm:continuous-ideal-weight} is framed in the Hilbert space setting of \eqref{eq:staf-problem-rkhs}. Thus, Theorem \ref{thm:continuous-ideal-weight} together with Theorem \ref{thm:feasibility} provides the StaF framework for RKHSs.

\begin{theorem}\label{thm:continuous-ideal-weight}Let $H$ be a RKHS over a set $X \subset \mathbb{R}^n$ with a strictly positive kernel $K:X\times X \to \mathbb{C}$ such that $K(\cdot,c) \in C^{m_0}(\mathbb{R}^n)$ for all $c \in X$. Suppose that $V \in H$. Let $C$ be an ordered collection of $M$ distinct centers, $C=(c_1,c_2,...,c_M) \in X^M$, with the associated ideal weights
\begin{equation}\label{eq:ideal-weight-center}W(C) = \argmin_{(a_i)_{i=1}^M \in \mathbb{C}^M} \left\| \sum_{i=1}^M a_i K(\cdot,c_i) - V(\cdot) \right\|_H.\end{equation} The function $W$ is $m_0$-times continuously differentiable with respect to each component of $C$.\end{theorem}

\begin{pf}The determination of $W(C)$ is equivalent to computing the projection of $V$ onto the space $Y=\vspan\{ K(\cdot, c_i) : i = 1,...,M \}$. To compute the projection, a Gram-Schmidt algorithm can be employed. The Gram-Schmidt algorithm is most easily expressed in its determinant form. Let $D_0=1$ and $D_m = \det \left( K(c_j,c_i) \right)_{i,j=1}^m$, then for $m=1,...,M$ the functions $$u_m(x) := \frac{1}{\sqrt{D_{m-1}D_m}} \det 
\begin{pmatrix}
K(c_1,c_1) & K(c_1,c_2) & \cdots & K(c_1,c_m)\\
K(c_2,c_1) & K(c_2,c_2) & \cdots & K(c_2,c_m)\\
\vdots & \vdots & \ddots & \vdots\\
K(c_{m-1},c_1) & K(c_{m-1},c_2) & \cdots & K(c_{m-1},c_m)\\
K(x, c_1) & K(x, c_2) & \cdots & K(x, c_m)
\end{pmatrix}$$
constitute an orthonormal basis for $Y$. Since $K$ is strictly positive definite, $D_m$ is positive for each $m$ and every $C$. The coefficient for each $K(x,c_l)$ with $l = 1,...,m$ in $u_m$ is a sum of products of the terms $K(c_i,c_j)$ for $i,j=1,...m$. Each such coefficient is $m_0$-times differentiable with respect to each $c_i$, $i=1,...,M$. When $\langle V, u_m \rangle$ is computed for the projection, the result is a linear combination of evaluations of $V$ at each of the centers. The function $V$ is $m_0$-times continuously differentiable, since $K$ is $m_0$-times differentiable \cite{Steinwart.Christmann2008}, therefore $\langle V, u_m \rangle$ is continuous with respect to the centers. Finally, each term in $W(C)$ is a linear combination of the coefficients determined by $u_m$ for $m=1,...,M$, and thus is $m_0$-times continuously differentiable with respect to each $c_i$ for $i=1,...,M$.
\qed\end{pf}

\section{The Gradient Chase Theorem}
\label{sec:gradient-chase}

As mentioned before, control theory problems involving function approximation and system identification are centered around the concept of weight update laws. Weight update laws are a collection of rules that the approximating weights must obey which lead to convergence to the ideal weights. In the case of the StaF approximation framework, the ideal weights are replaced with ideal weight functions. Theorem \ref{thm:continuous-ideal-weight} showed that if the moving centers of the StaF kernel functions are selected in such a way that the centers adjust smoothly with respect to the state $x$, then the ideal weight functions will also change smoothly with respect to $x$. Thus, in this context, weight update laws of the StaF approximation framework aim to achieve an estimation of the ideal weight function at the current state.

Theorem \ref{thm:gradient-chase} provides an example of such weight update laws that achieve a UUB result. The theorem takes advantage of perfect samples of a function in the RKHS $H$ corresponding to a real valued kernel function.

The proof of the theorem is similar to the standard proof for the convergence of the gradient descent algorithm for a quadratic programming problem \cite{Bertsekas99}. The contribution of the proof is in a modification, where the mean value theorem is used to produce an extra term which results in a UUB result, and the continuity of the largest and smallest eigenvalues of a Gram matrix are used to get a uniform bound in tandem with the Kantorovich inequality.

\begin{theorem}[Gradient Chase Theorem]\label{thm:gradient-chase}Let $H$ be a real valued RKHS over $\mathbb{R}^n$ with a continuously differentiable strictly positive definite kernel function $K:\mathbb{R}^n \times \mathbb{R}^n \to \mathbb{R}$. Let $V \in H$, $D \subset \mathbb{R}^n$ be a compact set, and $x:\mathbb{R}\to\mathbb{R}^n$ a state variable subject to the dynamical system $\dot x = q(x,t)$, where $q:\mathbb{R}^n \times \mathbb{R}_+ \to \mathbb{R}^n$ is a bounded locally Lipschitz continuous function. Further suppose that $x(t) \in D$ for all $t > 0$. Let $c:\mathbb{R}^n \to \mathbb{R}^M$, where for each $i=1,...,M$, $c_i(x)=x+d_i(x)$ where $d_i \in C^1(\mathbb{R}^n)$, and let $a \in \mathbb{R}^M$. Consider the function
$$F(a,c) = \left\| V - \sum_{i=1}^M a_i K(\cdot, c_i(x)) \right\|_H^2.$$
At each time instance $t >0$, there is a unique $W(t)$ for which $W(t) = \argmin_{a\in\mathbb{R}^M} F(a,c(x(t)))$. Given any $\epsilon >0$ and initial value $a^0$, there is a frequency $\tau > 0$, where if the gradient descent algorithm (with respect to $a$) is iterated at time steps $\Delta t < \tau^{-1}$, then $F(a^k,c^k) - F(w^k, c^k)$ will approach a neighborhood of radius $\epsilon$ as $k \to \infty$.\end{theorem}

\begin{pf}
Let $\bar \epsilon > 0$. By the Hilbert space structure of $H$: $$F(a,c) = \|V\|_H^2 -2V(c)^Ta + a^T K(c)a$$ where $V(c) = (V(c_1),...,V(c_M))^T$ and $K(c) = (K(c_i,c_j))_{i,j=1}^M$ is the symmetric strictly positive kernel matrix corresponding to $c$.
At each time iteration $t^k$, $k=0,1,2,...$, the corresponding centers and weights can be written as $c^k \in \mathbb{R}^{nM}$ and $a^k \in \mathbb{R}^M$, respectively. The ideal weights corresponding to $c^k$ will be denoted by $w^k$. It can be shown that $w^k = K(c^k)^{-1} V(c^k)$ and $F(w^k,c^k) = \|V\|_H^2 - V(c^k)^T K(c^k) V(c^k).$
Theorem \ref{thm:continuous-ideal-weight} ensures that the ideal weights change continuously with respect to the centers which remain in a compact set $\tilde D^M$, where $\tilde D = \{ x \in \mathbb{R}^M : \|x-D\| \le \max_{i=1,...,M} \left(\sup_{x \in D} |d_i(x)|\right) \}$, so the collection of ideal weights is bounded. Let $R > \bar \epsilon$ be large enough so that $N_R(0)$ contains both the initial value $a^0$ and the set of ideal weights. To facilitate the subsequent analysis, consider the constants:
\begin{align*}
R_0 &= \max_{x \in D, t > 0} |q(x,t)|&
R_1 &= \max_{a \in \overline{N_r(0)}, c \in \tilde D} |\nabla_a F(a,c)|\\
R_2 &= \max_{c \in \tilde D} |\nabla_c F(w(c),c)|&
R_3 &= \max_{c \in \tilde D} |\dot d_i(x(t)|\\
R_4 &= \max_{c \in \tilde D} \left\| \frac{d}{dc} w(c) \right\|
\end{align*}
where $\nabla_a$ is the gradient with respect to $a$, and let $\Delta t < \tau^{-1} :=\bar \epsilon \cdot (2(R_0+R_3)\cdot(R_1\cdot R_4\cdot (R_0+R_3)+ R_2+ 1))^{-1}$.
The proof aims to show that by using the gradient descent law for choosing $a^k$, the following inequality can be achieved:
$$\frac{F(a^{k+1},c^{k+1})-F(w^{k+1},c^{k+1})}{F(a^k,c^k)-F(w^k,c^k)} < \delta + \frac{\bar \epsilon}{F(a^k,c^k)-F(w^k,c^k)}$$
for some $0<\delta < 1$. Set \begin{equation}\label{eq:gradient-descent}a^{k+1} = a^k + \lambda g\end{equation} where $g=-\nabla_a F(a^k,c^k) = 2V(c^k)-2K(c^k)a^k$ and $\lambda$ is selected so that the quantity $F(a^k + \lambda g, c^k)$ is minimized. The $\lambda$ that minimizes this quantity is $\lambda = \left( \frac{g^Tg}{2g^TK(c^k)g}\right)$ which yields $F(a^{k+1},c^k) = F(a^k,c^k) - \frac{(g^Tg)^2}{4g^TK(c^k)g}$.
Since $F(a^{k+1},c^{k+1})$ is continuously differentiable in the second variable, we have $F(a^{k+1}, c^{k+1}) = F(a^{k+1},c^{k}) + \nabla_c F(a^{k+1},\eta)\cdot (c^{k+1}-c^k)$. Since $|\dot c(x(t))| < R_0 + R_3$, an application of the mean value theorem demonstrates that $\|c^{k+1} - c^k\| < (R_0+R_3)\Delta t$. Thus $$F(a^{k+1},c^{k+1}) = F(a^{k+1},c^k) + \epsilon_1(t^k),$$ where $|\epsilon_1(t^k)| < \bar \epsilon/2$ for all $k$. The quantity $F(w^{k+1},c^{k+1})$ is continuously differentiable in both variables. Thus, by the multi-variable chain rule and another application of the mean value theorem: $$F(w^{k+1},c^{k+1}) = F(w^{k},c^k)+\epsilon_2(t^k),$$ for $|\epsilon_2(t^k)| < \bar \epsilon/2$ for all $k$. Therefore, the following is established: 
\begin{align*}\frac{F(a^{k+1},c^{k+1}) - F(w^{k+1},c^{k+1})}{F(a^k,c^k)-F(w^k,c^k)} = \frac{F(a^{k+1},c^k)-F(w^{k},c^k) + (\epsilon_1(t^k)-\epsilon_2(t^k))}{F(a^k,c^k)-F(w^k,c^k)}\\
=1-\frac{(g^Tg)^2}{(g^TK(c^k)g)(g^TK(c^k)^{-1}g)} + \frac{\epsilon_1(t^k) - \epsilon_2(t^k)}{F(a^k,c^k)-F(w^k,c^k)}.\end{align*}

The Kantorovich inequality \cite{Bertsekas99} yields \begin{equation}\label{eq:kant-ineq}1-\frac{(g^Tg)^2}{(g^TK(c^k)g)(g^TK(c^k)^{-1}g)} \le \left(\frac{A_{c^k}/a_{c^k} - 1}{A_{c^k}/a_{c^k} + 1}\right)^2\end{equation} where $A_{c^k}$ is the largest eigenvalue of $K(c^k)$ and $a_{c^k}$ is the smallest eigenvalue of $K(c^k)$. The quantity on the right of \eqref{eq:kant-ineq} is continuous with respect to $A_{c^k}$ and $a_{c^k}$. In turn, $A_{c^k}$ and $a_{c^k}$ are continuous with respect to $K(c^k)$ (c.f. Exercise 4.1.6 \cite{Pedersen.Pedersen.ea1989}) which is continuous with respect to $c^k$. Therefore there is a largest value, $\delta$, that the right hand side of \eqref{eq:kant-ineq} obtains on the compact set $\tilde D$ and this value is less than $1$. Moreover, $\delta$ is independent of $\bar \epsilon$, so it may be declared that $\bar \epsilon = \epsilon(1-\delta)$. Finally,
$$\frac{F(a^{k+1},c^{k+1}) - F(w^{k+1},c^{k+1})}{F(a^k,c^k)-F(w^k,c^k)} \le \delta + \frac{(\epsilon_1(t^k) - \epsilon_2(t^k))}{F(a^k,c^k)-F(w^k,c^k)}.$$ Therefore, setting $e(k)=F(a^k,c^k)-F(w^k,c^k)$, it can be shown that $e(k+1) \le \delta e(k) + \epsilon(1-\delta)$ and the conclusion of the theorem follows.
\qed\end{pf}

\section{Simulation for the Gradient Chase Theorem}
\label{sec:gradient-chase-simulation}

To demonstrate the effectiveness of the Gradient Chase theorem, a simulation performed on a two-dimensional linear system is presented below. The system dynamics are given by
$$\begin{pmatrix}\dot x_1\\ \dot x_2\end{pmatrix} = \begin{pmatrix}0&1\\-1&0\end{pmatrix}\begin{pmatrix}x_1\\x_2\end{pmatrix},$$
which is the dynamical system corresponding to a circular trajectory. The function to be approximated is
$$V(x_1,x_2) = x_1^2 + 5x_2^2 + \tanh(x_1\cdot x_2),$$
and the kernel function to be used for function approximation are the exponential kernels, $K(x,y)=\exp\left(x^Ty\right)$. The centers are arranged in an equilateral triangle centered about the state. In particular, each center resides on a circle of radius $0.1$ centered at the state: $$c_i(x) = x + 0.1 \begin{pmatrix}\sin((i-1)2\pi/3)\\ \cos((i-1)2\pi/3)\end{pmatrix}$$ for $i=1,2,3$.

\begin{figure}[p] %H declares it must be HERE, h is a soft placement, and p gives the figure its own page
	\subfloat[][Trajectory of the state vector.]{
		\includegraphics[scale=.4]{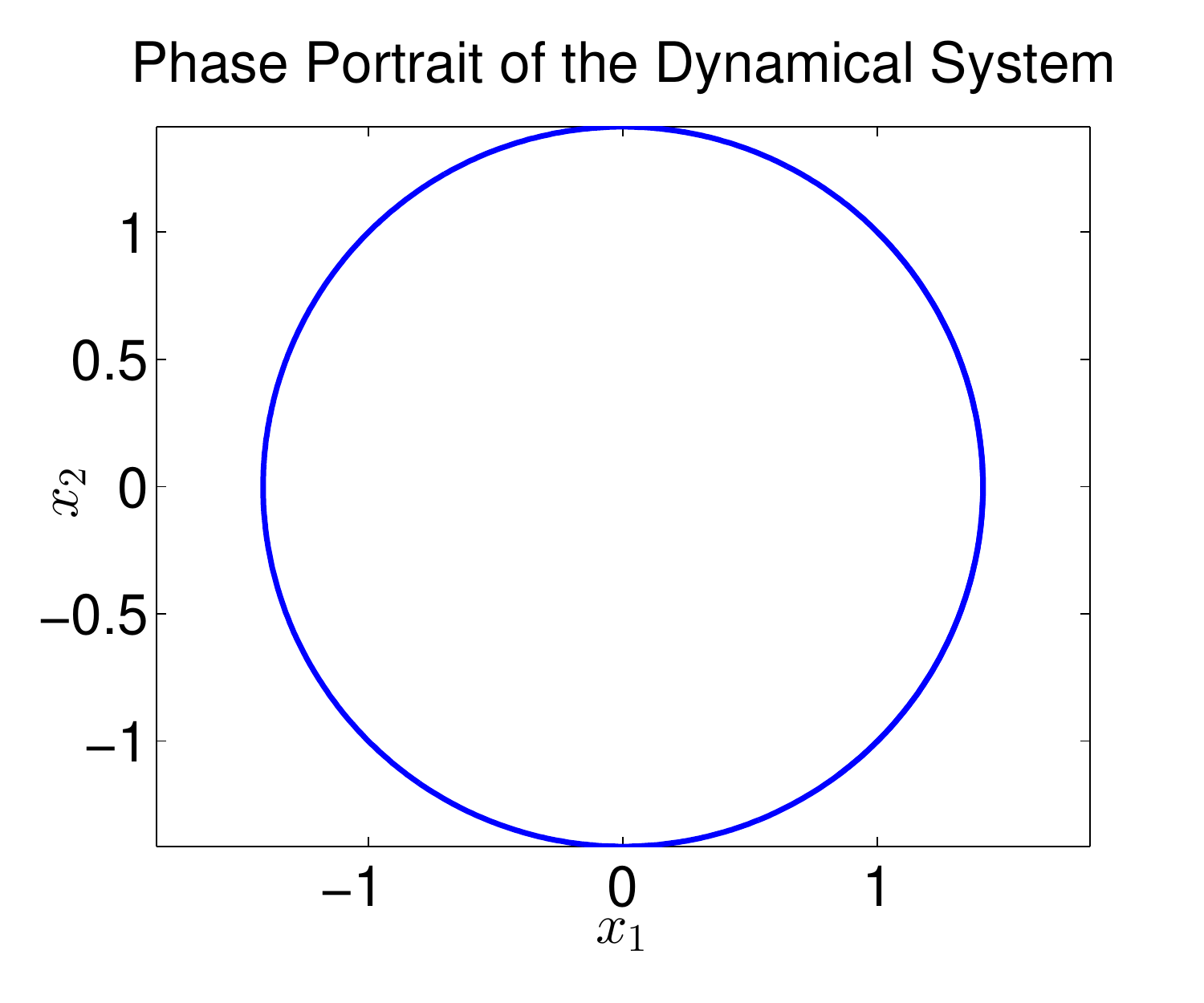}
		\label{fig:GD_StateTrajectory}
	}
	\subfloat[][Comparison of $V$ and the approximation $\hat V$.]{
		\includegraphics[scale=.4]{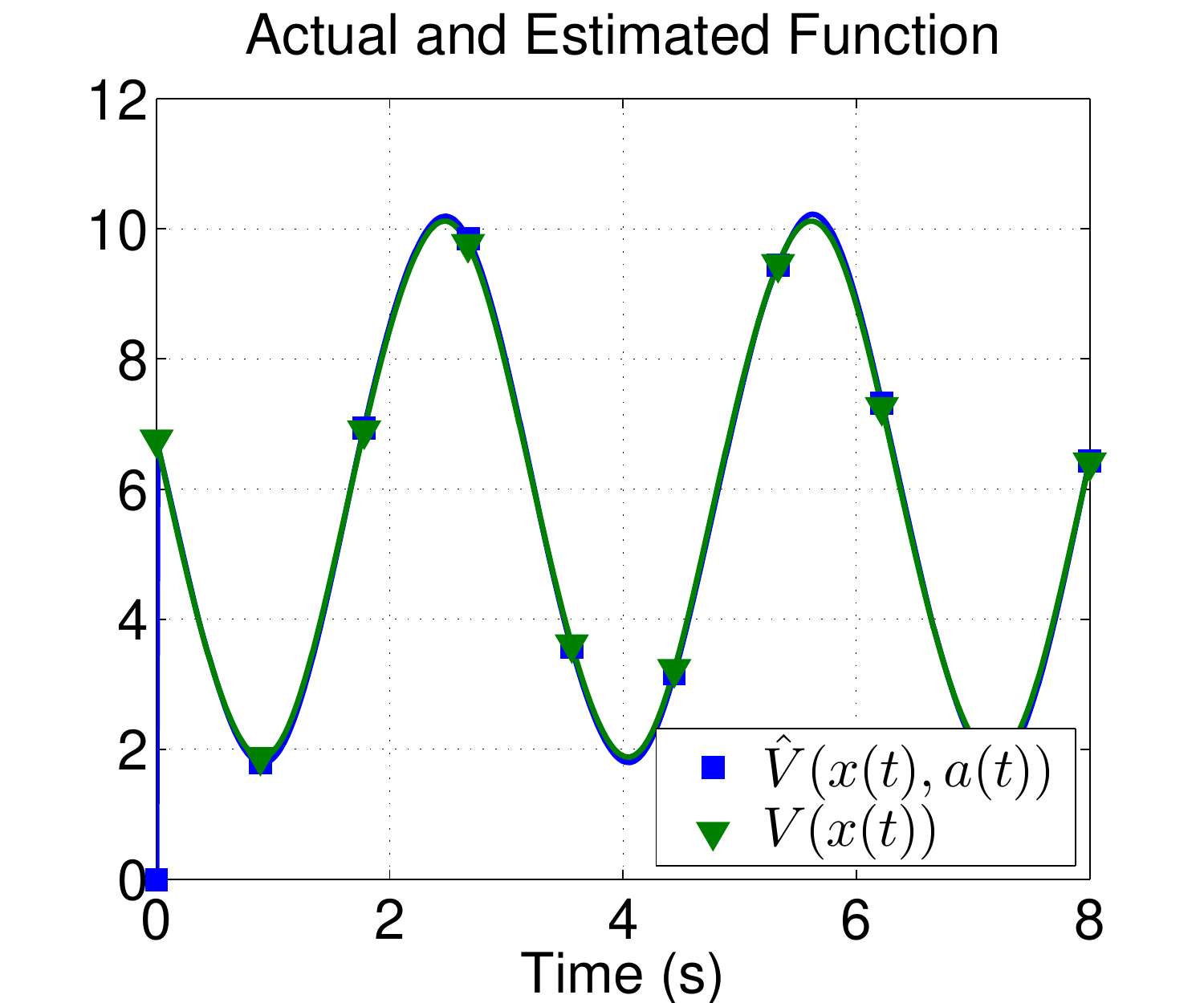}
		\label{fig:GD_Comparison}
	}
	
	\subfloat[][The values of the weight function estimates.]{
		\includegraphics[scale=.4]{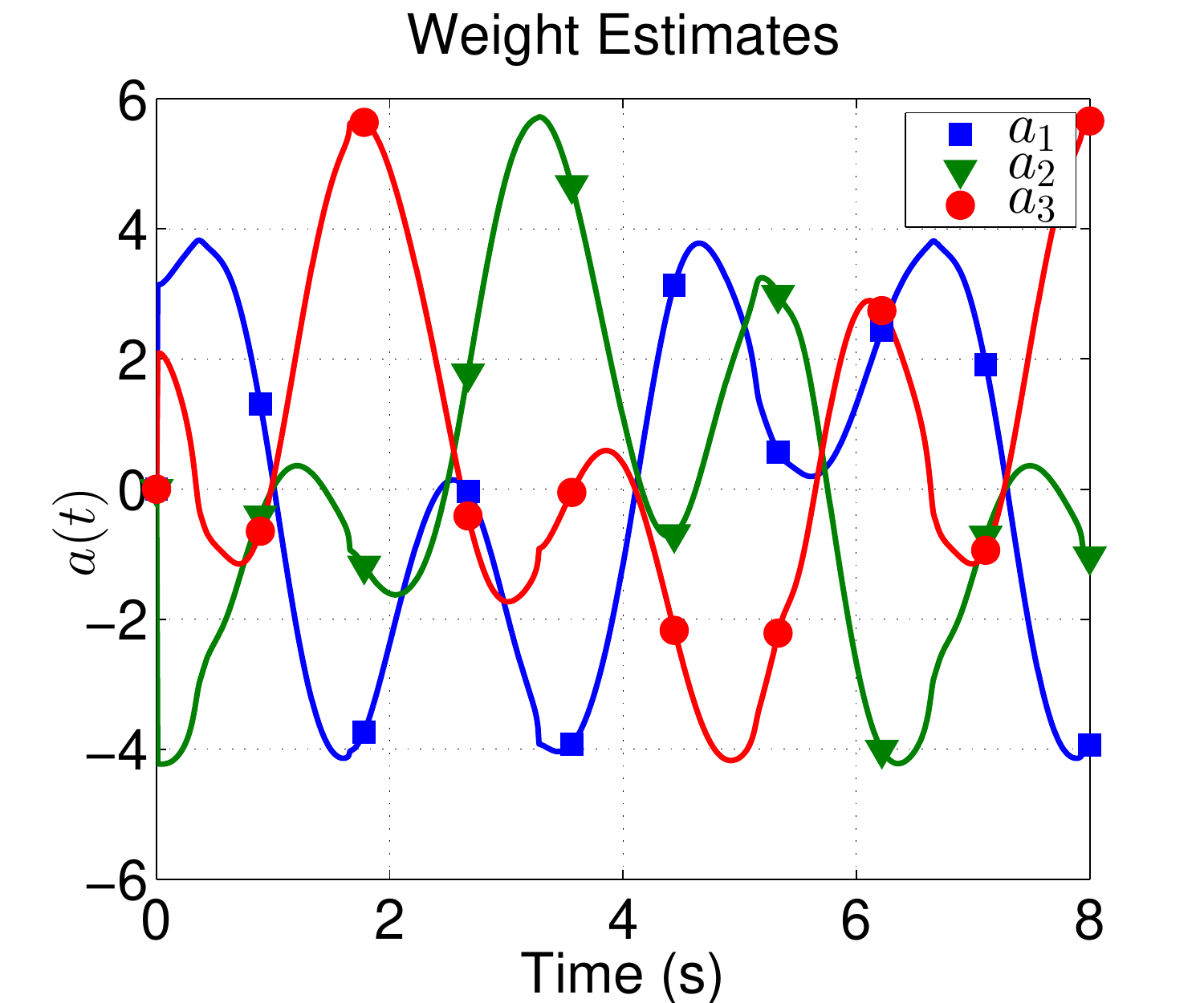}
		\label{fig:GD_Weights}
	}
	\subfloat[][Error committed by the approximation at the current state.]{
		\includegraphics[scale=.4]{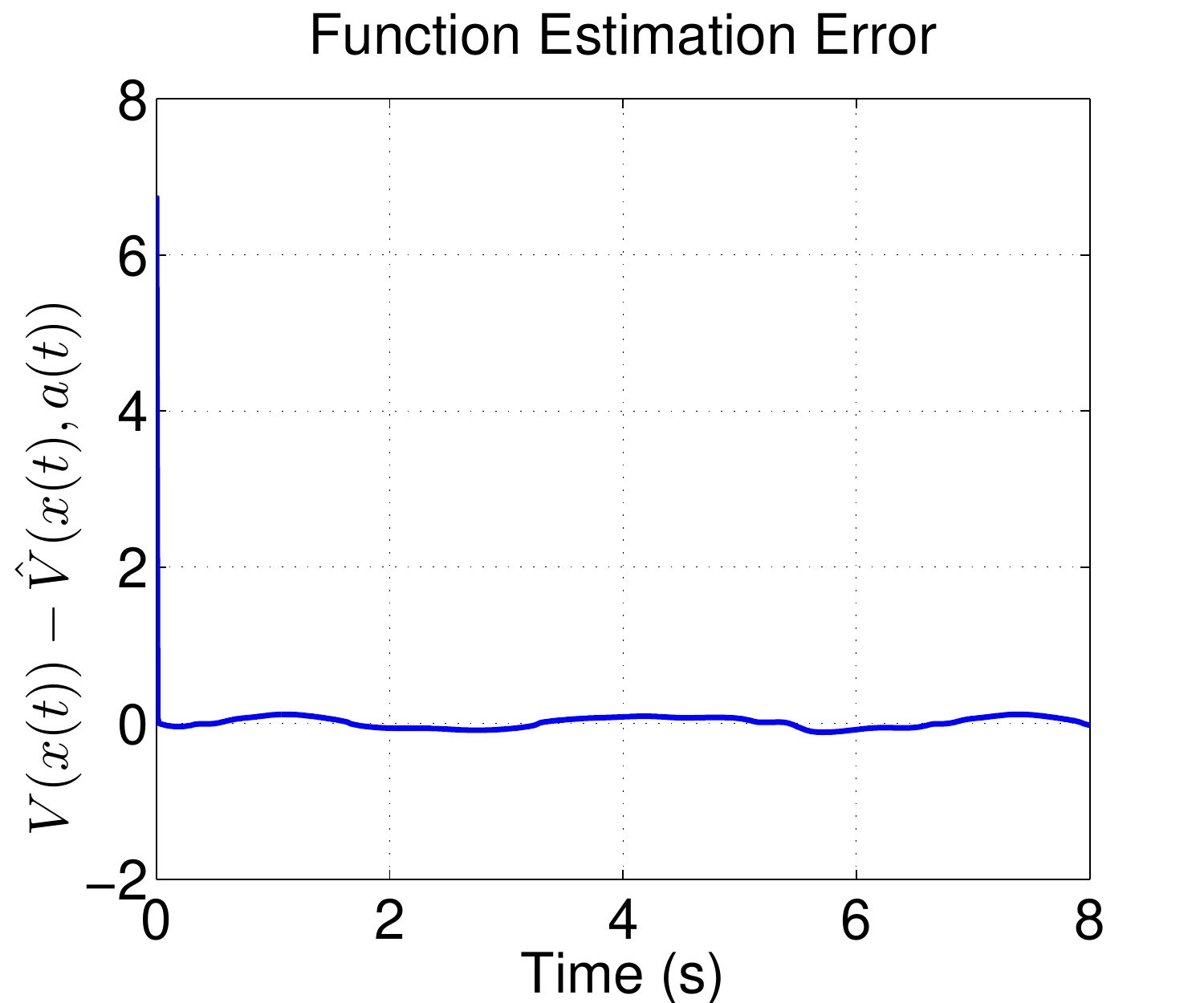}
		\label{fig:GD_Error}
	}
	\caption{Results of the numerical experiment demonstrating the Gradient Chase algorithm.}
	\label{fig:GD_Simulation}
\end{figure}

The initial values selected for the weights are $a^0 = [0\ 0\ 0]^T$. The gradient descent weight update law, given by \eqref{eq:gradient-descent}, are applied $10$ iterations per time-step and the time-steps incremented every $0.01$ seconds. Figure \ref{fig:GD_Simulation} presents the results of the simulation.

Figure \ref{fig:GD_Error} demonstrates that the function approximation error is regulated to a small neighborhood of zero as the Gradient Chase Theorem is implemented and validates the claim of the UUB result of Theorem \ref{thm:gradient-chase}. In Figure \ref{fig:GD_Weights}, approximations of the ideal weight function can be seen to be periodic as well as smooth. The smoothness of the ideal weight function itself is given in Theorem \ref{thm:continuous-ideal-weight}, and the periodicity of the approximation follows from the periodicity of the selected dynamical system, as illustrated in Figure \ref{fig:GD_StateTrajectory}. Figure \ref{fig:GD_Comparison} presents a comparison of $V$ evaluated at the current state to that of the approximation evaluated at the current state. Approximation of the function is maintained as the system state moves through its domain as anticipated.

\section{Application to Adaptive Dynamic Programming}\label{sec:ADP-Outline}

The application of approximation theory to the theory of optimal control arises through the approximation of the optimal value function, which is the solution to the Hamilton-Jacobi-Bellman (HJB) equation. Efficient methods for the approximation of the optimal value function are essential, since an increase in dimension can lead to a exponential increase in the number of required basis functions necessary to achieve an accurate approximation, the so called ``curse of dimensionality''.

The optimal value function corresponds to the infinite horizon optimal regulator problem, where the cost function $$J(x,u) = \int_0^\infty x^TQx + u^TRu\, dt$$ is to be minimized subject to the dynamics
\begin{equation}\label{eq:dynamics}\dot x(t) = f(x(t)) + g(x(t))u(t)\end{equation} where $x: \mathbb{R}_+ \to \mathbb{R}^n$, $u : \mathbb{R}_+ \to \mathbb{R}^m$, $Q \in \mathbb{R}^{n\times n}$, $R \in \mathbb{R}^{m\times m}$, with $Q$ and $R$ positive definite, $f: \mathbb{R}^n \to \mathbb{R}^n$, $g: \mathbb{R}^n \to \mathbb{R}^{n\times m}$. Moreover, $f$ and $g$ are assumed to be locally Lipschitz. The optimal value function is given by $$V(x) = \inf_{u \in \mathcal{U}} \int_0^\infty x^T Qx + u^T Ru\, dt$$ where $\mathcal{U}$ is the collection of admissible controllers. When the optimal value function is continuously differentiable and an optimal controller, $u^* \in \mathcal{U}$ exists, the optimal value function is the unique solution to the HJB equation \begin{equation}\label{eq:HJB}0 = x^TQx + u^{*T}Ru^* + \nabla V(x)(f(x) + g(x)u^*).\end{equation}

Once the optimal value function is determined, the optimal controller takes the form \begin{equation}\label{eq:optimal-controller}u^*(x(t)) = -\frac12 R^{-1} g(x)^T \nabla V(x(t))^T.\end{equation} In many applications, an approximation of the optimal controller is used real-time to yield autonomous behavior in a dynamic environment.

%Reinforcement learning (RL) is an approach that is frequently used for the determining online solutions of optimal control problems for systems with finite state and action spaces. In recent years, adaptive dynamic programming (ADP) has been successfully used to implement RL in deterministic autonomous control-affine systems (such as \eqref{eq:dynamics}) to solve optimal control problems via approximations of the optimal value function \cite{Al-Tamimi2008, Bhasin.Kamalapurkar.ea2013a, Dierks2009, Lewis.Vrabie2009, Mehta.Meyn2009, Padhi2006, Vamvoudakis2010, Zhang.Cui.ea2013, Zhang.Cui.ea2011, Zhang.Liu.ea2013}. ADP techniques aim to approximate the value function through a parameterization by collections of basis functions and through the real-time tuning of approximation weights. ADP differs from approximate dynamic programming (also refered to as ADP in the literature), in that it uses online adaptive control techniques to achieve an approximation and to attain the stability of a dynamical system.

For some problems, such as the linear quadratic regulator (LQR) problem, the optimal value function takes a particular form which simplifies the choice of basis functions. In the case of LQR, the optimal value function is of the form $\sum_{i,j=1}^n w_{i,j} x_j x_i$ (c.f. \cite{Kirk2004,Liberzon2012}), so basis functions of the form $\sigma_{i,j} = x_j x_i$ will provide an accurate estimation of the optimal value function provided the weights, $w_{i,j} \in \mathbb{R}$, are tuned properly. However, in most cases, the form of the optimal value function is unknown, and generic basis functions have been proposed to parameterize the problem.

Adaptive dynamic programming (ADP) replaces $V$ with a parametrization, $\hat V(x, W_c) = \sum_{i=1}^M w_{i,c} \sigma_i(x)$, with $W_c = (w_{1,c},..., w_{M,c}) \in \mathbb{R}^M$, and $u^*$ with a parametrization $\hat u(x, W_a) = -\frac12 R^{-1} g(x)^T \nabla_x V(x, W_a)^T$ where $W_a \in \mathbb{R}^M$. The actor and critic weights, $W_a$ and $W_c$ respectively, are tuned to minimize the residual Bellman error (BE), $$\delta(x, W_a, W_c) = x^TQ x + \hat u(x,W_a)^T R u(x,W_a) + \nabla_x \hat V(x,W_c) \left( f(x) + g(x) \hat u(x,W_a)\right),$$ over all $x$ in some compact set $D$ in real-time. The BE is used to motivate weight update laws for $W_a$ and $W_c$ to achieve a real-time minimization.

Throughout the ADP literature, many basis functions have been proposed for real-time (approximate) optimal control. However, in practice, it is difficult to select weight update laws that guarantee stability by achieving a good approximation of the ideal weights, especially for a system with a modest embedded processor. In the majority of cases, actual implementation of ADP is executed using only polynomial basis functions, and the StaF method enables a broader class of functions to be used for approximate optimal control of a dynamical system.

In this setting, the StaF problem becomes $$\limsup_{t\to \infty} \sup_{x \in \overline{N_r(x)}} |\delta(x, W_a(t), W_c(t))| < \epsilon.$$

\ref{sec:adp-simulation} provides more information concerning the application of the StaF method to ADP by presenting the results of a companion paper \cite{Kamalapurkar.Rosenfeld.ea2015}.

%Often, kernel functions from reproducing kernel Hilbert spaces (RKHSs) are used as generic basis functions, and the approximation problem is solved over a (preferably large) compact domain of $\mathbb{R}^n$ \cite{Christmann.Christmann.ea2010,Micchelli.Micchelli.ea2007,Park.Park.ea1991}. An essential property of RKHSs is given a collection of basis functions in the Hilbert space, there is a unique set of weights that minimize the error in the Hilbert space norm, the so called ideal weights \cite{Folland1999}. The modal choice of kernel as a basis function is that of the Gaussian radial basis function (RBF) given by $K(x,y) = \exp(-\|x-y\|^2/\mu)$ where $x,y \in \mathbb{R}^n$ and $\mu > 0$ \cite{Park1991,Steinwart.Christmann2008}. For the approximation of a function over a large compact domain, a large number of basis functions is required, which leads to an intractable computational problem for online control applications.

\section{Conclusion} A new StaF kernel method is introduced in this paper for the purpose of function approximation. The development in this paper establishes that by using the StaF method a local approximation of a function can be maintained in real-time as a state moves through a compact domain. Heuristically, much fewer kernel functions are required in comparison to more traditional function approximation schemes, since the approximation is maintained in a smaller region. For the exponential kernels, a new theorem in this paper establishes that an explicit bound on the number of kernel functions required can calculated. Two applications of this methodology were presented. In Section \ref{sec:gradient-chase}, a ``gradient chase'' algorithm was developed. There it was seen that a function may be well approximated provided that the algorithm was applied with a high enough frequency. Simulations results provided in Section \ref{sec:gradient-chase-simulation} demonstrated the performance of the gradient chase algorithm, and an application to ADP is provided in Section \ref{sec:ADP-Outline} and the Appendix for an infinite horizon optimal regulation problem.

The strength of the StaF methodology is the reduction of the computational requirements for real-time implementation of a function approximation, through the reduction in the number of basis functions. As demonstrated in \ref{sec:adp-simulation}, where only three basis functions were required to achieve a stabilizing approximate optimal controller for a $2$-dimensional system. However, since the StaF method aims at maintaining an accurate local approximation of the value function only in a local neighborhood of the current system state, the StaF kernel method lacks memory, in the sense that the information about the ideal weights over a region of interest is lost when the state leaves the region of interest. Thus, unlike existing techniques, the StaF method generates an approximation that is valid only in a local region. A memory-based modification to the StaF kernel method that retains and reuses past information for creating a global approximation is the subject of future research.

\appendix
\section{Applications to Adaptive Dynamic Programming}
\label{sec:adp-simulation}

\begin{figure}[p]
	\subfloat[][Trajectory of the state vector.]{
		\includegraphics[scale=.4]{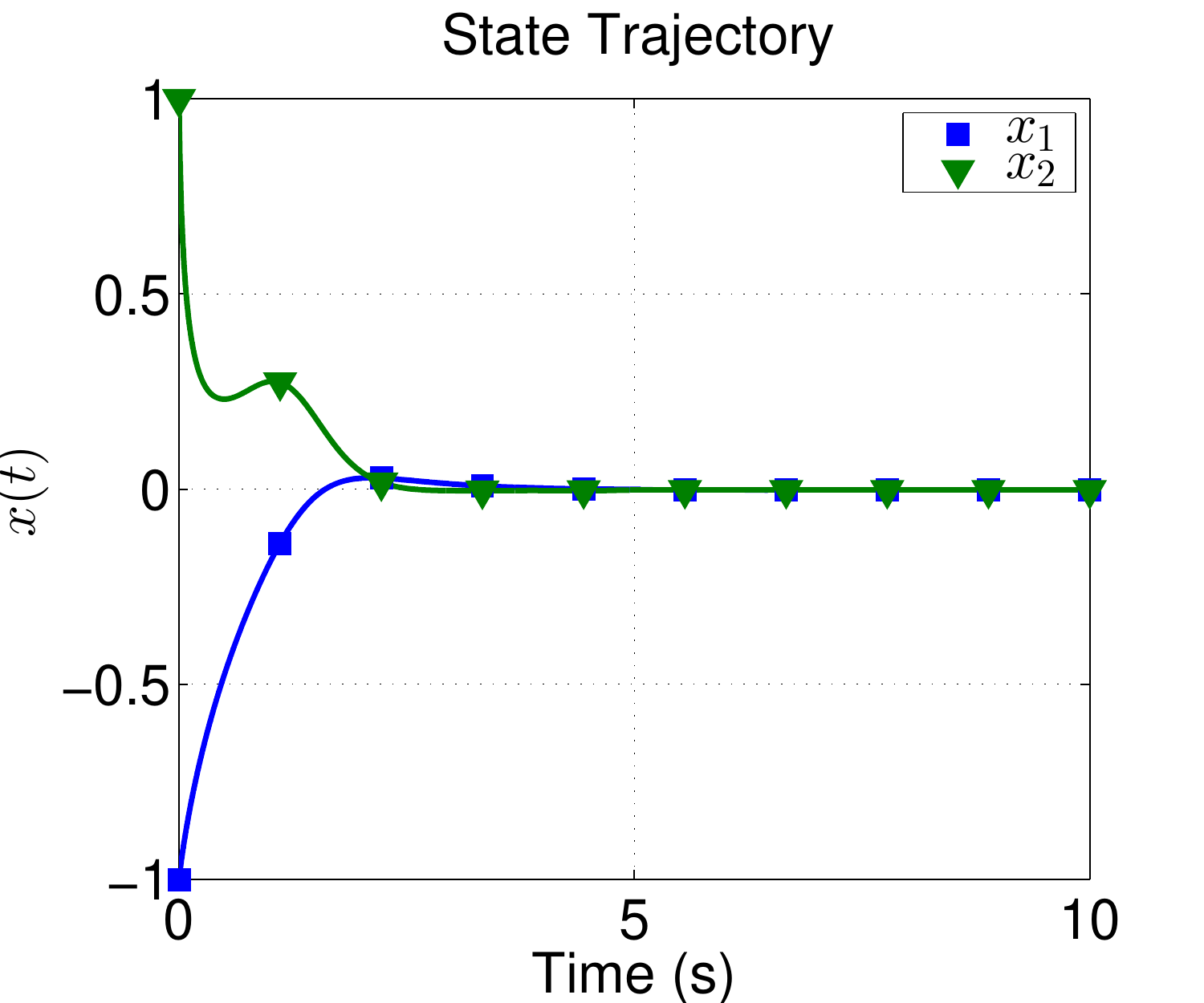}
		\label{fig:ADP-State}
	}
	\subfloat[][Trajectory of the actor weights.]{
		\includegraphics[scale=.4]{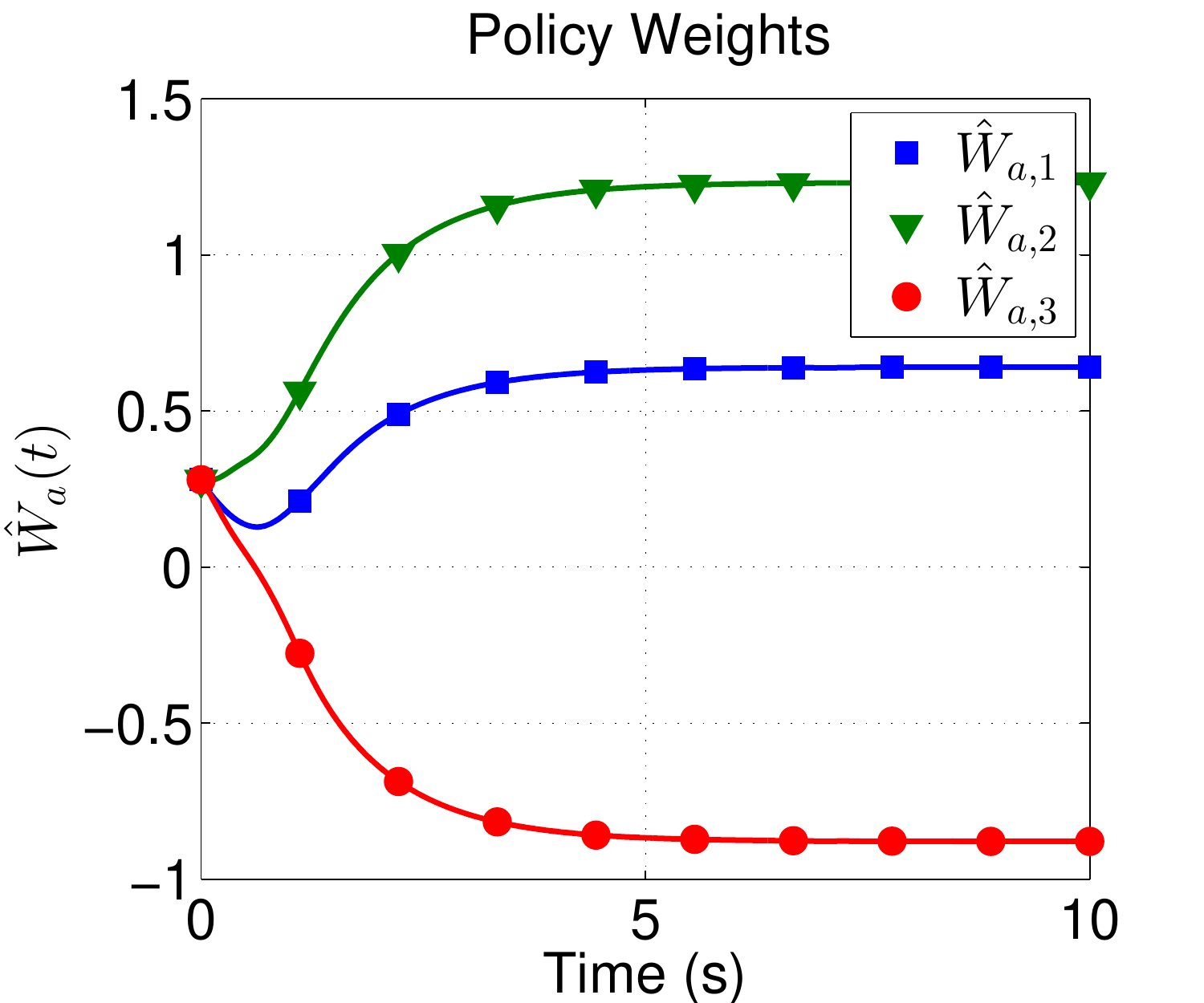}
		\label{fig:ADP-ActorWeights}
	}
	
	\subfloat[][The values of the critic weights.]{
		\includegraphics[scale=.4]{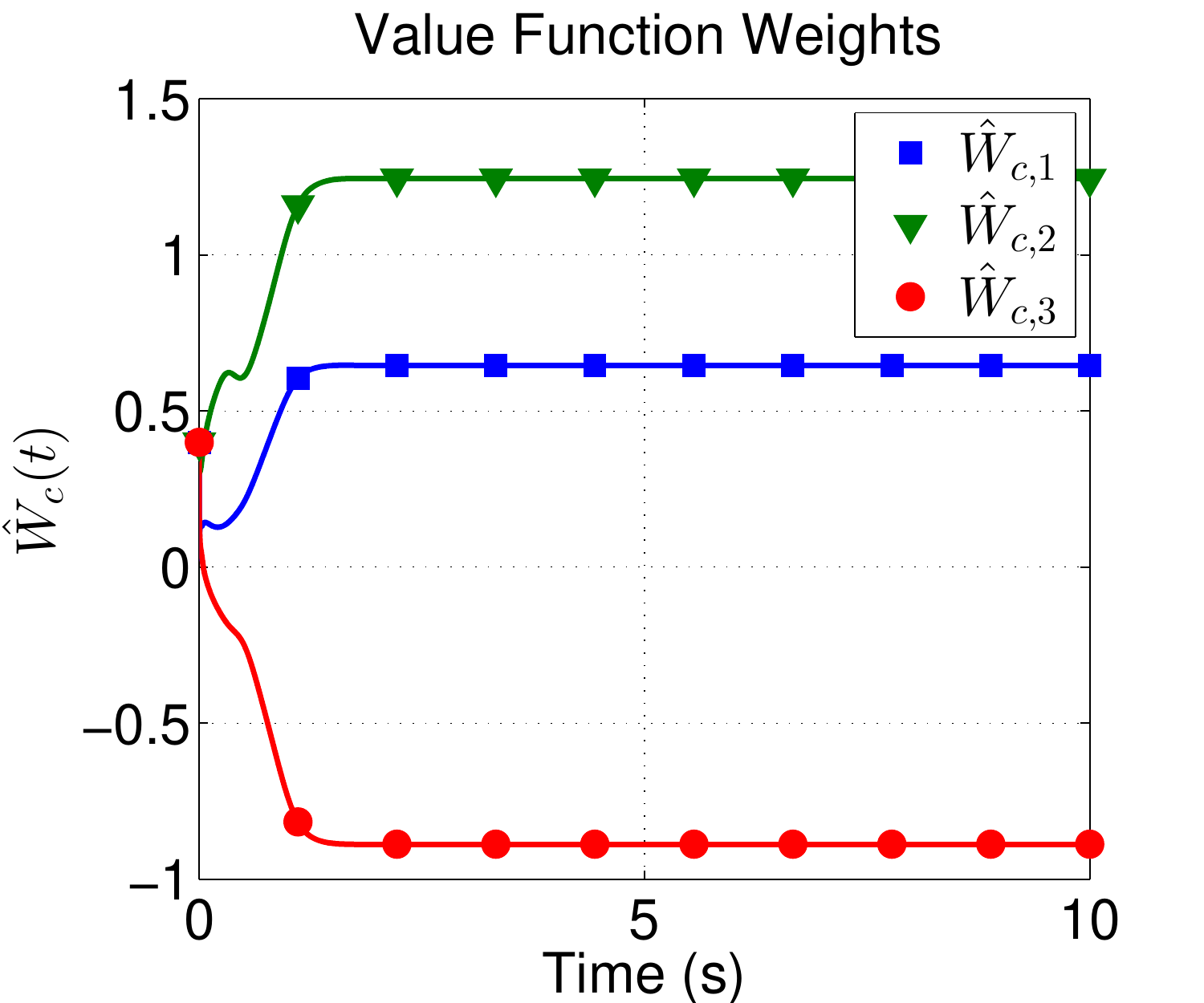}
		\label{fig:ADP-CriticWeights}
	}
	\subfloat[][Error committed by approximate policy.]{
		\includegraphics[scale=.4]{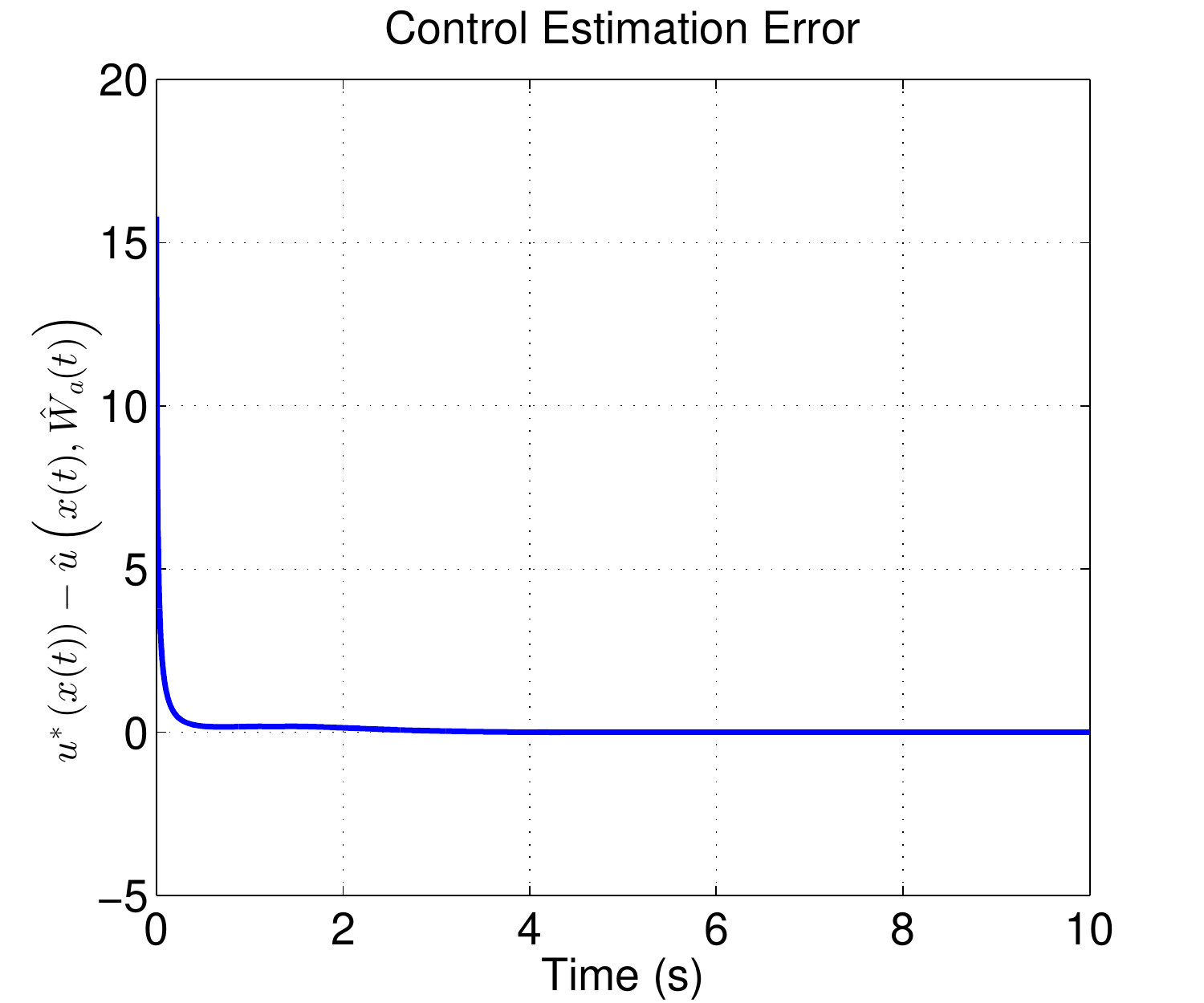}
		\label{fig:ADP-Opt-Control-Error}
	}
	
	\subfloat[][Error of the estimation of the value function at the current state.]{
		\includegraphics[scale=.4]{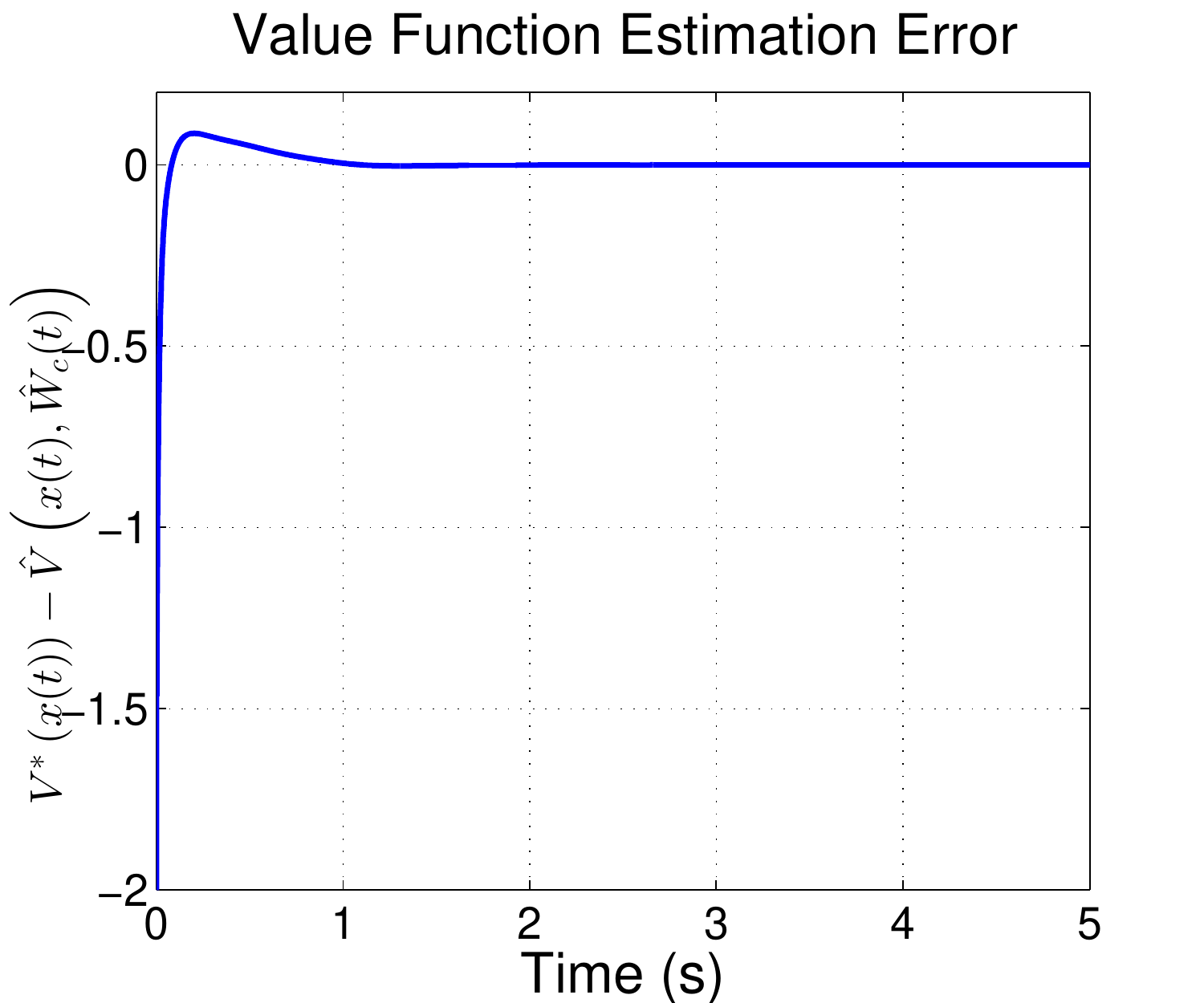}
		\label{fig:ADP-Value-Error}
	}
	\caption{Results of the numerical experiment demonstrating the convergence for the StaF ADP method. \cite{Kamalapurkar.Rosenfeld.ea2015}}
	\label{fig:ADP-Simulation}
\end{figure}

To demonstrate the effectiveness of the StaF technique in the context of optimal control, the simulation results of a companion paper are presented here. The details of the analysis are contained in \cite{Kamalapurkar.Rosenfeld.ea2015}.
The dynamical system in question is of the form $\dot x = f(x)+g(x)u$ where $x=(x_1,x_2)^T \in \mathbb{R}^2$,
\begin{equation}\label{eq:optimal-regulator}f(x) = \begin{pmatrix} x_2-x_1\\ -\frac12 x_1 - \frac12 x_2(\cos(2x_1)+2)^2\end{pmatrix}\text{, and }g(x) = \begin{pmatrix}0\\ \cos(2x_1)+2\end{pmatrix}.\end{equation}
Associated with this dynamical system is the cost functional
\begin{equation}\label{eq:cost}J(x,u)=\int_0^\infty (x^T(\tau)x(\tau) + u(\tau)^2)\,d\tau\end{equation}
In the infinite horizon regulation problem, the goal is to determine an optimal control law $u^*:\mathbb{R}^2 \to \mathbb{R}$ (assuming an optimal control law exists) that satisfies $$u^*(x_0) =\argmin_{u\in\mathcal{U}} \int_0^\infty (x^T(\tau)x(\tau)+u(\tau)^2)\,d\tau$$ where $\mathcal{U}$ is the collection of admissible controllers and $x(0)=x_0$ inside the integrand. The optimal value function is given by $$V(x_0) = \min_{u\in\mathcal{U}} \int_0^\infty (x^T(\tau)x(\tau)+u(\tau)^2)\,d\tau$$ when such a minimum exists, and the optimal value function satisfies the HJB equation \eqref{eq:HJB}. If $V^*$ satisfies the HJB equation and is also continuously differentiable, then it is the unique solution to \eqref{eq:HJB}. Furthermore, $u^*$ can be determined from $V^*$ by $u^*(x) = \frac12 g^T(x) \nabla V^*(x)$.

In most cases, the optimal value function cannot be determined analytically, and approximate solutions are used instead. However, for the system presented in the section, the optimal value function is known. In particular, for the infinite horizon optimal regulator problem with dynamics given by \eqref{eq:optimal-regulator} with cost functional \eqref{eq:cost}, the optimal value function is given by $V^*(x) = \frac12 x_1^2 + x_2^2$ and the associated optimal control law is given by $u^*(x) = -(\cos(2x_1)+2)x_2.$ More details can be found in \cite{Vamvoudakis2010}.

In this example, the infinite horizon optimal regulator problem is solved in real-time. The function $V^*$ is approximated by a function of the form $$\hat V(x,\hat W_c) =\sum_{i=1}^3 \hat W_{c,i} \left(\exp(x^T c_i(x)) - 1 \right)$$ where $\hat W_c \in \mathbb{R}^3$ are weights to be adjusted in real-time, and $c_i(x) = x + d_i(x)$ where $$d_i(x) = 0.7 \left( \frac{x^Tx + 0.01}{1+x^Tx} \right) \left[ \cos\left( \frac{2\pi}{3} i + \frac{\pi}{2} \right), \sin\left( \frac{2\pi}{3} i + \frac{\pi}{2} \right) \right]^T$$ for $i=1,2,3$. The approximation of the optimal control law is given by $$\hat u(x,\hat W_a) = -\frac12 g^T(x) \nabla_x \hat V(x,\hat W_a)$$ where $\hat W_a \in \mathbb{R}^3$ are weights to be adjusted in real-time.
In the framework of ADP, the functions $V^*$ and $u^*$ are replaced by their approximations $\hat V$ and $\hat u$, respectively, in the HJB equation, yielding a residual nonzero error, called the Bellman error (BE). The goal is to minimize the BE by adjustments of the weights, $\hat W_a$ and $\hat W_c$. If the BE is identically zero after the adjustment of the weights, then the optimal value function and the approximation of the optimal value function coincide. For nonzero BE, the BE is used as a heuristic measure of the distance between $\hat V$ and $V^*$, as well as the distance between $\hat u$ and $u^*$. The weight update laws and subsequent convergence analysis can be found in \cite{Kamalapurkar.Rosenfeld.ea2015}.

The results of the numerical experiment are presented in Figure \ref{fig:ADP-Simulation}. Figure \ref{fig:ADP-State} indicates that the state is regulated to the origin when using the ADP algorithm combined with the StaF methodology. Figure \ref{fig:ADP-ActorWeights} shows that the weight vector $\hat W_a$ converged as well. In typical StaF implementations, the weights are not expected to converge. However, since the optimal control problem is a regulator problem, the state and the centers ultimately occupy a fixed neighborhood of the origin, and the weights converge to the ideal weights corresponding to a small neighborhood of the origin.

When the weights converge, it is expected that $\hat W_a$ and $\hat W_c$ converge to the same values. The convergence is demonstrated by comparing Figure \ref{fig:ADP-ActorWeights} and Figure \ref{fig:ADP-CriticWeights}. The approximate controller and the optimal controller converge as well, as shown in Figure \ref{fig:ADP-Opt-Control-Error}, and the value function estimation error, given in Figure \ref{fig:ADP-Value-Error}, vanishes rapidly.
\section*{Acknowledgments}

This research is supported in part by NSF award number 1509516 and Office of Naval Research Grant N00014-13-1-0151. Any opinions, findings and conclusions or recommendations expressed in this material are those of the author(s) and do not necessarily reflect the views of the sponsoring agencies.

\bibliographystyle{elsarticle-num}
\bibliography{master,ncr,encr}

\begin{thebibliography}{10}
\expandafter\ifx\csname url\endcsname\relax
  \def\url#1{\texttt{#1}}\fi
\expandafter\ifx\csname urlprefix\endcsname\relax\def\urlprefix{URL }\fi
\expandafter\ifx\csname href\endcsname\relax
  \def\href#1#2{#2} \def\path#1{#1}\fi

\bibitem{Rosenfeld.Kamalapurkar.ea2015}
J.~A. Rosenfeld, R.~Kamalapurkar, W.~E. Dixon, State following ({S}ta{F})
  kernel functions for function approximation {P}art {I}: Theory and
  motivation, in: Proc. Am. Control Conf., 2015, pp. 1217--1222.

\bibitem{Christmann.Christmann.ea2010}
A.~Christmann, I.~Steinwart, Univeral kernels on non-standard input spaces, in:
  Advances in Nueral Information Processing, 2010, pp. 406--414.

\bibitem{Micchelli.Micchelli.ea2007}
C.~A. Micchelli, Y.~Xu, H.~Zhang, Universal kernels, J. Mach. Learn. Res. 7
  (2006) 2651--2667.

\bibitem{Park.Park.ea1991}
J.~Park, I.~W. Sandberg,
  \href{http://dx.doi.org/10.1162/neco.1991.3.2.246}{Universal approximation
  using radial-basis-function networks}, Neural Comput. 3~(2) (1991) 246--257.
\newblock \href {http://dx.doi.org/10.1162/neco.1991.3.2.246}
  {\path{doi:10.1162/neco.1991.3.2.246}}.
\newline\urlprefix\url{http://dx.doi.org/10.1162/neco.1991.3.2.246}

\bibitem{C.Prudhomme.C.Prudhomme.ea2001}
C.~Prud'homme, D.~V. Rovas, K.~Veroy, L.~Machiels, Y.~Maday, A.~T. Patera,
  G.~Turinici, Reliable real-time solution of parametrized partial differential
  equations: Reduced-basis output bound methods, AMSE Journal of Fluids
  Engineering 124(1) (2001) 70--80.

\bibitem{Balmes.Balmes.ea1996}
E.~Balmes,
  \href{http://www.sciencedirect.com/science/article/pii/S0888327096900278}{Parametric
  families of reduced finite element models. theory and applications},
  Mechanical Systems and Signal Processing 10~(4) (1996) 381 -- 394.
\newblock \href {http://dx.doi.org/http://dx.doi.org/10.1006/mssp.1996.0027}
  {\path{doi:http://dx.doi.org/10.1006/mssp.1996.0027}}.
\newline\urlprefix\url{http://www.sciencedirect.com/science/article/pii/S0888327096900278}

\bibitem{AhmedK.Noor.AhmedK.Noor.ea1980}
J.~M.~P. Ahmed K.~Noor, Reduced basis technique for nonlinear analysis of
  structures, AIAA Journal 18~(2) (1980) 455--462.

\bibitem{Al-Tamimi2008}
A.~Al-Tamimi, F.~L. Lewis, M.~Abu-Khalaf, Discrete-time nonlinear {HJB}
  solution using approximate dynamic programming: Convergence proof, IEEE
  Trans. Syst. Man Cybern. Part B Cybern. 38 (2008) 943--949.

\bibitem{Bhasin.Kamalapurkar.ea2013a}
S.~Bhasin, R.~Kamalapurkar, M.~Johnson, K.~G. Vamvoudakis, F.~L. Lewis, W.~E.
  Dixon, A novel actor-critic-identifier architecture for approximate optimal
  control of uncertain nonlinear systems, Automatica 49~(1) (2013) 89--92.

\bibitem{Dierks2009}
T.~Dierks, B.~Thumati, S.~Jagannathan, {Optimal control of unknown affine
  nonlinear discrete-time systems using offline-trained neural networks with
  proof of convergence}, Neural Netw. 22~(5-6) (2009) 851--860.

\bibitem{Lewis.Vrabie2009}
F.~L. Lewis, D.~Vrabie, Reinforcement learning and adaptive dynamic programming
  for feedback control, IEEE Circuits Syst. Mag. 9~(3) (2009) 32--50.

\bibitem{Mehta.Meyn2009}
P.~Mehta, S.~Meyn, Q-learning and pontryagin's minimum principle, in: Proc.
  IEEE Conf. Decis. Control, 2009, pp. 3598 --3605.

\bibitem{Padhi2006}
R.~Padhi, N.~Unnikrishnan, X.~Wang, S.~Balakrishnan, {A single network adaptive
  critic (SNAC) architecture for optimal control synthesis for a class of
  nonlinear systems}, Neural Netw. 19~(10) (2006) 1648--1660.

\bibitem{Vamvoudakis2010}
K.~Vamvoudakis, F.~Lewis, {Online actor-critic algorithm to solve the
  continuous-time infinite horizon optimal control problem}, Automatica 46~(5)
  (2010) 878--888.

\bibitem{Zhang.Cui.ea2013}
H.~Zhang, L.~Cui, Y.~Luo, Near-optimal control for nonzero-sum differential
  games of continuous-time nonlinear systems using single-network adp, IEEE
  Trans. Cybern. 43~(1) (2013) 206--216.

\bibitem{Zhang.Cui.ea2011}
H.~Zhang, L.~Cui, X.~Zhang, Y.~Luo, Data-driven robust approximate optimal
  tracking control for unknown general nonlinear systems using adaptive dynamic
  programming method, IEEE Trans. Neural Netw. 22~(12) (2011) 2226--2236.

\bibitem{Zhang.Liu.ea2013}
H.~Zhang, D.~Liu, Y.~Luo, D.~Wang, Adaptive Dynamic Programming for Control
  Algorithms and Stability, Communications and Control Engineering,
  Springer-Verlag, London, 2013.

\bibitem{DeVore.DeVore.ea1998}
R.~A. DeVore, Nonlinear approximation, Acta Numerica 7 (1998) 51--150.

\bibitem{Gaggero.Gaggero.ea2013}
M.~Gaggero, G.~Gnecco, M.~Sanguineti, Dynamic programming and value-function
  approximation in sequential decision problems: Error analysis and numerical
  results, Journal of Optimization Theory and Applications 156.

\bibitem{Gaggero.Gaggero.ea2014}
M.~Gaggero, G.~Gnecco, M.~Sanguineti, Approximate dynamic programming for
  stochastic n-stage optimization with application to optimal consumption under
  uncertainty, Computational Optimization and Applications 58~(1) (2014)
  31--85.

\bibitem{Zoppoli.Zoppoli.ea2002}
R.~Zoppoli, M.~Sanguineti, T.~Parisini, Approximating networks and extended
  ritz method for the solution of functional optimization problems, Journal of
  Optimization Theory and Applications 112~(2) (2002) 403--440.
\newblock \href {http://dx.doi.org/10.1023/A:1013662124879}
  {\path{doi:10.1023/A:1013662124879}}.

\bibitem{Kamalapurkar.Rosenfeld.ea2015}
R.~Kamalapurkar, J.~A. Rosenfeld, W.~E. Dixon, State following ({S}ta{F})
  kernel functions for function approximation {P}art {II}: Adaptive dynamic
  programming, in: Proc. Am. Control Conf., 2015, pp. 521--526.

\bibitem{Zhu.Zhu.ea2012}
K.~Zhu, \href{http://dx.doi.org/10.1007/978-1-4419-8801-0}{Analysis on {F}ock
  spaces}, Vol. 263 of Graduate Texts in Mathematics, Springer, New York, 2012.
\newblock \href {http://dx.doi.org/10.1007/978-1-4419-8801-0}
  {\path{doi:10.1007/978-1-4419-8801-0}}.
\newline\urlprefix\url{http://dx.doi.org/10.1007/978-1-4419-8801-0}

\bibitem{Pinkus.Pinkus.ea2004}
A.~Pinkus, Strictly positive definite functions on a real inner product space,
  Adv. in Comput. Math. 20 (2004) 263--271.

\bibitem{Steinwart.Christmann2008}
I.~Steinwart, A.~Christmann, Support vector machines, Information Science and
  Statistics, Springer, New York, 2008.

\bibitem{Beylkin.Beylkin.ea2005}
G.~Beylkin, L.~Monzon,
  \href{http://www.sciencedirect.com/science/article/pii/S106352030500014X}{On
  approximation of functions by exponential sums}, Applied and Computational
  Harmonic Analysis 19~(1) (2005) 17 -- 48.
\newblock \href
  {http://dx.doi.org/http://dx.doi.org/10.1016/j.acha.2005.01.003}
  {\path{doi:http://dx.doi.org/10.1016/j.acha.2005.01.003}}.
\newline\urlprefix\url{http://www.sciencedirect.com/science/article/pii/S106352030500014X}

\bibitem{Bertsekas99}
D.~P. Bertsekas, Nonlinear Programming, Athena Scientific, Belmont, MA, 1999.

\bibitem{Pedersen.Pedersen.ea1989}
G.~K. Pedersen, \href{http://dx.doi.org/10.1007/978-1-4612-1007-8}{Analysis
  now}, Vol. 118 of Graduate Texts in Mathematics, Springer-Verlag, New York,
  1989.
\newblock \href {http://dx.doi.org/10.1007/978-1-4612-1007-8}
  {\path{doi:10.1007/978-1-4612-1007-8}}.
\newline\urlprefix\url{http://dx.doi.org/10.1007/978-1-4612-1007-8}

\bibitem{Kirk2004}
D.~Kirk, Optimal Control Theory: An Introduction, Dover, 2004.

\bibitem{Liberzon2012}
D.~Liberzon, Calculus of variations and optimal control theory: a concise
  introduction, Princeton University Press, 2012.

\end{thebibliography}
\end{document}